\documentclass[12pt,reqno]{amsart}
\usepackage{graphicx}
\usepackage{tabularx}
\usepackage{hyperref}
 \newtheorem{thm}{Theorem}[section]
 \newtheorem{cor}[thm]{Corollary}
 \newtheorem{lem}[thm]{Lemma}
 
 \theoremstyle{definition}
 
 \theoremstyle{remark}
 
 \numberwithin{equation}{section}

\begin{document}
\begin{center}
{\Large\textsc{\textbf{A Note on a result due to Ankeny and Rivlin\footnote{This is a preprint of a paper whose final and definite form is published open access in Applied Mathematics E-Notes. See http://www.math.nthu.edu.tw/~amen/ for the final version.}}}}
\vskip0.5in

\textbf{Eze R. Nwaeze}\\[0pt]
Department of Mathematics, Tuskegee University\\
Tuskegee, AL 36088, USA\\
\verb"enwaeze@mytu.tuskegee.edu"\\
\end{center}
\qquad \\

\begin{quotation}{\bf Abstract.} ~~Let $p(z)=a_0+a_1z+a_2z^2+a_3z^3+\cdots+a_nz^n$ be a polynomial of degree $n$ having no zeros in the unit disk. ~Then it is well known that for $R\geq 1,$
$\displaystyle{\max_{|z|=R}|p(z)|}\leq \Big(\dfrac{R^n+1}{2}\Big)\displaystyle{\max_{|z|=1}|p(z)|}.$ In this paper, we consider polynomials with gaps, having all its zeros on the circle $S(0, K):=\{z: |z|=K\}, ~0<K\le 1,$~ and estimate the value of $\Big(\dfrac{{\max_{|z|=R}|p(z)|}}{{\max_{|z|=1}|p(z)|}}\Big)^s$ for any positive integer $s.$ \\


\end{quotation}

\section{Introduction}

Let $p(z)=\displaystyle{\sum_{j=0}^{n}a_{j}z^{j}}$ be a polynomial of degree $n.$ We will denote
$$M(p, r):=\displaystyle{\max_{|z|=r}|p(z)|}\,, ~r>0,$$
$$||p||:= \displaystyle{\max_{|z|=1}|p(z)|}\,,$$ and
$$D(0, K):=\{z: |z|<K\},~ K> 0.$$
 Bernstein observed the following result, which in fact is a simple consequence of the maximum modulus principle  (see \cite[p.~137]{Polya}). This inequality is also known as the Bernstein's inequality.

\begin{thm}\label{thm1.1}
Let $p(z)=\displaystyle{\sum_{j=0}^{n}a_{j}z^{j}}$ be a polynomial of degree $n.$ Then for $R\geq 1,$
\begin{equation}
M(p,R)\leq R^n||p||.
\end{equation}
Equality holds for $p(z)=\alpha z^n,$ $\alpha$ being a complex number.
\end{thm}

For polynomial of degree $n$ not vanishing in the interior of the unit circle, Ankeny and Rivlin \cite{Ankeny} proved the following result.

\begin{thm}[\bf Ankeny and Rivlin \cite{Ankeny} ]\label{thm1.2}
Let $p(z)=\displaystyle{\sum_{j=0}^{n}a_{j}z^{j}}\neq 0$ in $D(0, 1).$ Then for $R\geq 1,$
\begin{equation}\label{eqn1}
M(p,R)\leq \Big(\dfrac{R^n+1}{2}\Big)||p||.
\end{equation}
Here equality holds for $p(z)=\dfrac{\alpha+\beta z^n}{2},$ where $|\alpha|=|\beta|=1.$
\end{thm}

In 2005, Gardner, Govil and Musukula \cite{Gardner2} proved the following generalization and sharpening of Theorem \ref{thm1.2}.

\begin{thm}\label{thm2.5}
Let $p(z)=a_0+\displaystyle{\sum_{j=t}^{n}a_{j}z^{j}}, ~1\leq t\leq n,$  be a polynomial of degree $n$ and $p(z)\neq 0$ in $D(0, K), ~K\geq 1.$ Then for $R\geq
1,$

\begin{eqnarray}\label{eqn9}
M(p,R)\leq \Big(\dfrac{R^n+s_0}{1+s_0}\Big)||p||-\Big(\dfrac{R^n-1}{1+s_0}\Big)m
- \dfrac{n}{1+s_0}\Bigg[\dfrac{(||p||-m)^2-(1+s_0)^{2}|a_n|^2}{(||p||-m)}\Bigg] \\ \nonumber
\times\Bigg\{ \dfrac{(R-1)(||p||-m)}{(||p||-m)+(1+s_0)|a_n|}-\ln\Bigg[1+\dfrac{(R-1)(||p||-m)}{(||p||-m)+(1+s_0)|a_n|}\Bigg]\Bigg\},\nonumber
\end{eqnarray}
where $m=\displaystyle{\min_{|z|=K}|p(z)|},$ and
$$s_0=K^{t+1}\dfrac{\frac{t}{n}\cdot\frac{|a_t|}{|a_0|-m}K^{t-1}+1}{\frac{t}{n}\cdot\frac{|a_t|}{|a_0|-m}K^{t+1}+1}.$$

\end{thm}

Several research monographs have been written on this subject of inequalities (see for example Govil and Mohapatra \cite{Govil6}, Milovanovi\'c, Mitrinovi\'c and Rassias \cite{MMR}, Rahman and Schmeisser \cite{Rahman}, and recent article of Govil and Nwaeze \cite{GovilNwaeze2}).\\

While trying to obtain an inequality analogous to (\ref{eqn1}) for polynomials not vanishing in $D(0, K), K\leq 1,$ Dewan and Ahuja \cite{Dewan2} were able to prove this only for polynomials having all the zeros on the circle $S(0, K):=\{z: |z|=K\}, ~0<K\le 1.$

\begin{thm}\label{thm4.1}
Let $p(z)=\displaystyle{\sum_{j=0}^{n}a_{j}z^{j}}$  be a polynomial of degree $n$ having all its zeros on $S(0, K), \; K\leq 1.$ Then for $R\geq 1$ and for every
positive integer $s,$

\begin{equation}
\{M(p, R)\}^s\leq \Bigg[\dfrac{K^{n-1}(1+K)+(R^{ns}-1)}{K^{n-1}+K^n}\Bigg]\{M(p, 1)\}^s.
\end{equation}
\end{thm}

For $s=1$, the Theorem \ref{thm4.1} reduces to

\begin{cor}
Let $p(z)=\displaystyle{\sum_{j=0}^{n}a_{j}z^{j}}$  be a polynomial of degree $n$ having all its zeros on $S(0, K), \; K\leq 1.$ Then for $R\geq 1,$
\begin{equation}\label{eqn25}
M(p, R)\leq \Bigg[\dfrac{K^{n-1}(1+K)+(R^{n}-1)}{K^{n-1}+K^n}\Bigg]M(p, 1).
\end{equation}
\end{cor}

In same spirit, we prove the following results

\section{Main Results}

\begin{thm}\label{thm0.1}

Let $p(z)=z^m\Bigg[a_{n-m}z^{n-m}+\displaystyle{\sum_{j=\mu}^{n-m}a_{n-m-j}z^{n-m-j}}    \Bigg],$ $1\leq \mu\leq n-m,$ $0\leq m\leq n-1,$ be a polynomial of degree $n,$ having $m-fold$ zeros at origin and remaining $n-m$ zeros on $S(0, K),$ $K\leq 1.$ Then for $R\geq 1$ and every positive  integer $s,$
\begin{equation}
[M(p,R)]^s\leq L(\mu; K, m, n,s)[M(p,1)]^s,
\end{equation}
where\\ $$L(\mu; K, m, n,s)=\dfrac{n(K^{n-m-2\mu+1}+K^{n-m-\mu+1})+(R^{ns}-1)[n+mK^{n-m-2\mu+1}+mK^{n-m-\mu+1}-m]}{n(K^{n-m-2\mu+1}+K^{n-m-\mu+1})}.$$
\end{thm}
For $m=0,$ we have
\begin{cor}\label{cor0.2}
Let $p(z)=a_{n}z^{n}+\displaystyle{\sum_{j=\mu}^{n}a_{n-j}z^{n-j}},$ $1\leq \mu\leq n,$  be a polynomial of degree $n,$ having all zeros on $|z|=K,$ $K\leq 1.$ Then for $R\geq 1$ and every positive integer $s,$
\begin{equation}
[M(p,R)]^s\leq L(\mu; K, n,s)[M(p,1)]^s, 
\end{equation}
where \\$$L(\mu; K, n,s)=\dfrac{K^{n-\mu}(K^{1-\mu}+K)+(R^{ns}-1)}{K^{n-2\mu+1}+K^{n-\mu+1}}.$$
\end{cor}
If we set $\mu=1$ into Corollary \ref{cor0.2}, we get the following result of Dewan and Ahuja \cite{Dewan2}.
\begin{cor}\label{cor0.3}
Let $p(z)=\displaystyle{\sum_{j=0}^{n}a_{j}z^{j}},$ be a polynomial of degree $n,$ having all zeros on $|z|=K,$ $K\leq 1.$ Then for $R\geq 1$ and every positive integer $s,$
\begin{equation}
[M(p,R)]^s\leq L(1; K, n,s)[M(p,1)]^s, 
\end{equation}
where\\
$$L(1; K, n,s)=\dfrac{K^{n-1}(1+K)+(R^{ns}-1)}{K^{n-1}+K^{n}}.$$
\end{cor}

\section{Lemmas}
For the proof Theorem \ref{thm0.1} we need the following lemmas. The first lemma is due to Kumar and Lal \cite{SL}.
\begin{lem}\label{lem0.6}
Let $p(z)=z^m\Bigg[a_{n-m}z^{n-m}+\displaystyle{\sum_{j=\mu}^{n-m}a_{n-m-j}z^{n-m-j}}    \Bigg],$ $1\leq \mu\leq n-m,$~$0\leq m\leq n-1,$ be a polynomial of degree $n,$ having $m-fold$ zeros at origin and remaining $n-m$ zeros on $|z|=K,$ $K\leq 1.$
\begin{equation}
\displaystyle{\max_{|z|=1}|p'(z)|}\leq \dfrac{n+m(K^{n-m-2\mu+1}+K^{n-m-\mu+1}-1)}{K^{n-m-2\mu+1}+K^{n-m-\mu+1}}\displaystyle{\max_{|z|=1}|p(z)|}.
\end{equation}
\end{lem}
The next lemma is the Bernstein inequality given in Theorem \ref{thm1.1}.
\begin{lem}\label{lem0.7}
Let $p(z)$ be a polynomial of degree $n.$ Then for $R\geq 1,$
\begin{equation}
M(p,R)\leq R^{n}M(p,1).
\end{equation}
\end{lem}

\section{Proof}
\begin{proof}[\bf Proof of Theorem \ref{thm0.1}]

By Lemma \ref{lem0.6}, we have
$$\displaystyle{\max_{|z|=1}|p'(z)|}\leq \dfrac{n+m(K^{n-m-2\mu+1}+K^{n-m-\mu+1}-1)}{K^{n-m-2\mu+1}+K^{n-m-\mu+1}}\displaystyle{\max_{|z|=1}|p(z)|}.$$

Applying Lemma \ref{lem0.7} to the polynomial $p'(z)$ which is  of degree $n-1,$ it follows that for all $R\geq 1$ and $\theta \in [0,2\pi),$

\begin{align*}
|p'(Re^{i\theta})| &\leq \displaystyle{\max_{|z|=R}|p'(z)|}\\
&\leq R^{n-1}\displaystyle{\max_{|z|=1}|p'(z)|}\\
&\leq  R^{n-1}\Bigg[\dfrac{n+m(K^{n-m-2\mu+1}+K^{n-m-\mu+1}-1)}{K^{n-m-2\mu+1}+K^{n-m-\mu+1}}\Bigg]\displaystyle{\max_{|z|=1}|p(z)|}.
\end{align*}
So for each $\theta \in [0,2\pi)$ and $R\geq 1,$ we obtain

\begin{align*}
\big[p(Re^{i\theta})\big]^s - \big[p(e^{i\theta})\big]^s &= \int_{1}^{R}\dfrac{d\big[p(te^{i\theta})\big]^s}{dt}dt\\
&=\int_{1}^{R}s\big[p(te^{i\theta})\big]^{s-1}p'(te^{i\theta})e^{i\theta}dt.
\end{align*}

This implies that
$$\big|p(Re^{i\theta})\big|^s \leq \big|p(e^{i\theta})\big|^s + s\int_{1}^{R}\big|p(te^{i\theta})\big|^{s-1}\big|p'(te^{i\theta})\big|dt.$$
So,
\begin{align*}
\big[M(p,R)\big]^s &\leq \big[M(p,1)\big]^s + s\int_{1}^{R}\big[t^nM(p,1)\big]^{s-1}\big|p'(te^{i\theta})\big|dt\\
&\leq \big[M(p,1)\big]^s + s\int_{1}^{R}t^{ns-n}\big[M(p,1)\big]^{s-1}t^{n-1}\dfrac{n+m(K^{n-m-2\mu+1}+K^{n-m-\mu+1}-1)}{K^{n-m-2\mu+1}+K^{n-m-\mu+1}}M(p,1)dt\\
&=\big[M(p,1)\big]^s + s\Bigg[\dfrac{n+m(K^{n-m-2\mu+1}+K^{n-m-\mu+1}-1)}{K^{n-m-2\mu+1}+K^{n-m-\mu+1}}\Bigg]\big[M(p,1)\big]^s\int_{1}^{R}t^{ns-1}dt\\
&= \big[M(p,1)\big]^s + \big[M(p,1)\big]^{s}\Bigg[\dfrac{n+m(K^{n-m-2\mu+1}+K^{n-m-\mu+1}-1)}{K^{n-m-2\mu+1}+K^{n-m-\mu+1}}\Bigg]s\dfrac{R^{ns}-1}{ns}\\
&=  \big[M(p,1)\big]^s \Bigg[1+ \dfrac{\Big[n+m(K^{n-m-2\mu+1}+K^{n-m-\mu+1}-1) \Big]\big(R^{ns}-1\big)}{n\big(K^{n-m-2\mu+1}+K^{n-m-\mu+1}\big)}\Bigg].
\end{align*}
This yields
{\small $$\big[M(p,R)\big]^s\leq \big[M(p,1)\big]^s \Bigg[\dfrac{n\big(K^{n-m-2\mu+1}+K^{n-m-\mu+1}\big)+\Big[n+m(K^{n-m-2\mu+1}+K^{n-m-\mu+1}-1) \Big]\big(R^{ns}-1\big)}{n\big(K^{n-m-2\mu+1}+K^{n-m-\mu+1}\big)}\Bigg].$$}
This completes the proof.

\end{proof}

\textbf{Acknowledgment.} Many thanks to the anonymous referee for his/her valuable comments.

\end{document}